\documentclass{amsart}

\usepackage{amsmath,amssymb,amsthm,latexsym}

\newtheorem{theorem}{Theorem}
\newcommand{\bt}{\begin{theorem}}
\newcommand{\et}{\end{theorem}}
\newtheorem{lemma}{Lemma}
\newcommand{\bl}{\begin{lemma}}
\newcommand{\el}{\end{lemma}}
\newtheorem{corollary}{Corollary}
\newcommand{\bc}{\begin{corollary}}
\newcommand{\ec}{\end{corollary}}
\newcommand{\bconj}{\begin{conjecture}}
\newcommand{\econj}{\end{conjecture}}
\newtheorem{problem}{Problem}
\newcommand{\bprob}{\begin{problem}}
\newcommand{\eprob}{\end{problem}}
\newcommand{\beq}{\begin{equation}}
\newcommand{\eeq}{\end{equation}}
\newcommand{\benum}{\begin{enumerate}}
\newcommand{\eenum}{\end{enumerate}}
\newcommand{\N}{\ensuremath{ \mathbf N }}
\newcommand{\Z}{\ensuremath{\mathbf Z}}
\newcommand{\Q}{\ensuremath{\mathbf Q}}
\newcommand{\R}{\ensuremath{\mathbf R}}

\newcommand{\mcj}{\ensuremath{ \mathcal J}}

\newcommand{\mbe}{\ensuremath{ \mathbf e}}

\newcommand{\mbo}{\ensuremath{ \mathbf 0}}

\newcommand{\mbx}{\ensuremath{ \mathbf x}}
\newcommand{\mby}{\ensuremath{ \mathbf y}}
\newcommand{\mbz}{\ensuremath{ \mathbf z}}

\newcommand{\lex}{\ensuremath{\preceq_{\text{lex}}}}

\newcommand{\bmat}{\left(\begin{matrix}}
\newcommand{\emat}{\end{matrix}\right)}
\newcommand{\bsmallmat}{\left(\begin{smallmatrix}}
\newcommand{\esmallmat}{\end{smallmatrix}\right)}

\DeclareMathOperator{\qqand}{\qquad\text{and}\qquad}

\DeclareMathOperator{\vectorxn}{\left( \begin{matrix} x_1 \\ \vdots \\ x_n \end{matrix}\right)}
\DeclareMathOperator{\vectoryn}{\left( \begin{matrix} y_1 \\ \vdots \\ y_n \end{matrix}\right)}
\DeclareMathOperator{\vectorzn}{\left( \begin{matrix} z_1 \\ \vdots \\ z_n \end{matrix}\right)}

\title{Addition theorems in partially ordered  groups}
\author{Melvyn B. Nathanson}
\address{Department of Mathematics\\Lehman College (CUNY)\\Bronx, NY 10468}
\email{melvyn.nathanson@lehman.cuny.edu}

\date{\today}


 \thanks{Supported in part by  PSC-CUNY Research Award Program grant 66197-00 54.}

\begin{document}

\begin{abstract}
Shnirel'man's inequality and Shnirel'man's basis theorem are fundamental 
results about sums of sets of positive integers in additive number theory.   
It is proved that these results are inherently order-theoretic 
and extend to partially ordered abelian and nonabelian groups.  
One abelian application is an addition theorem 
for sums of sets of $n$-dimensional lattice points. 
\end{abstract}

\maketitle

\section{Shnirel'man density and Shnirel'man's inequality}
Let $\N = \{1,2,3,\ldots\}$ be the set of positive integers and 
$\N_0 = \{0,1,2,3,\ldots\}$ the set of nonnegative integers. 

Let $A$ be a  set  of positive integers.  
The \emph{counting function} of the set $A$  is 
\[
A(n) = \sum_{\substack{a \in A \\ 1 \leq a \leq n}} 1.
\]
The \emph{Shnirel'man density} of the set $A$ is 
\[  
\sigma(A) = \inf \left\{ \frac{A(n)}{n} : n \in \N \right\}.
\]
It follows that $\sigma(A) = 1$ if and only if $A = \N$. 

Shnirel'man density has some expected  properties, such as the inequality 
\[
0 \leq \sigma(A) \leq 1
\]
for every set $A$.  Also, $\sigma(A) = 1/2$ if $A$ is the set of odd 
integers.  Shnirel'man density also has some unusual properties.
For example, $\sigma(A) = 0$ if $A$ is the set of even integers.  
More generally, $\sigma(A) = 0$  if $1 \notin A$, and so 
\[
 \sigma(A) > 0 \quad  \text{ implies } \quad 1 \in A. 
\] 

For sets $A$ and $B$ of positive integers, we define the \emph{sumset} 
\[
A+B =  \{a+b:a \in A\cup \{0\} \text{ and } b \in B \cup \{0\} \}.
\]
The sum of $h$ sets $A_1, \ldots, A_h$ of positive integers is 
\[
A_1+\cdots + A_h = \left\{ a_1+\cdots + a_h: 
a_i \in A_i \cup \{0\} \text{ for all } i \in \{1,\ldots,h\} \right\}.
\] 
In particular, if $A_i = A$ for all  $i \in \{1,\ldots,h\}$, we have the \emph{$h$-fold sumset} 
\[
hA = \{a_1+\cdots + a_h: a_i \in A\cup \{0\} \text{ for all } i \in \{1,\ldots, h\} \}.
\]

The set $A$ is a \emph{basis of order $h$}  for the  nonnegative integers 
if $hA= \N_0$.  
The set $A$ is a \emph{basis} for the  nonnegative integers 
if $A$ is a  basis of order $h$ for some $h \geq 1$.  
Note that $A$ is a  basis of order $h$ if and only if $\sigma(hA) = 1$.
This is one reason why Shnirel'man density is a good density for many additive problems.

\bt                               \label{Ndensity:theorem:Shnirel-pigeonhole}
Let $A$ and $B$ be sets of positive integers.  
If $\sigma(A)+\sigma(B) \geq 1$, then $A+B = \N_0$. 
\et

\begin{proof}
The proof is a simple pigeonhole principle argument. 
Let $\sigma(A) = \alpha$ and $\sigma(B) = \beta$. 
We have  $0+0 \in A+B$. 
Because $\alpha > 0$ or $\beta > 0$, we have $1 \in A\cup B$ 
and $1 = 1+0 \in A+B$.  

Let $n \geq 2$.  If $n \in A \cup B$, then $n = n+0 \in A+B$. 
Suppose that $n \notin A\cup B$.  Let $A' = A \cap \{1,\ldots, n-1\}$ 
and $B' = B \cap \{1,\ldots, n-1\}$.  Then  
\[
n - B' = \{n-b:b \in B' \} \subseteq  \{1,\ldots, n-1\}. 
\]   
Because $n \notin A\cup B$, we have 
\[
|A'| =  A(n-1)  =  A(n)  \geq \alpha n 
\]
and 
\[
|n-B'| = |B'| = B(n-1) = B(n) \geq \beta n
\] 
and so 
\[
 |A'| + |n-B'|  \geq (\alpha + \beta) n \geq n.
 \]
Because $A'$ and $B'$ are subsets of $\{1,\ldots, n-1\}$, 
the sets $A'$ and $n-B'$ are not disjoint and so there exist  
$a \in A' \subseteq A$ and $b \in B' \subseteq B$ such that $a = n-b$.  
 Thus, $n = a+b \in A+B$.  
This completes the proof. 
\end{proof}

\bt[Shnirel'man's inequality]                             \label{Ndensity:theorem:Shnirel-ineq-classical}
Let $A$ and $B$ be sets of positive integers.    
If $\sigma(A) = \alpha$ and $\sigma(B) = \beta$, then 
\beq                                 \label{Ndensity:Shnirel-ineq-classical}
\sigma(A+B) \geq \alpha + \beta - \alpha\beta
\eeq
\et

\begin{proof}
This  is a special case of Theorem~\ref{Ndensity:theorem:LatticeShnirelman} below.
\end{proof}

Note that $0 \leq \alpha, \beta \leq 1$ implies 
\[
0 \leq  \alpha + \beta - \alpha\beta =  \alpha + \beta (1 - \alpha) \leq 1.
\]

The inequality 
\[
\sigma(A+B) \geq \alpha + \beta - \alpha\beta
\]
is equivalent to 
\[
1-\sigma(A+B) \leq (1-\alpha)(1-\beta). 
\]
It follows by induction on $h$ that if $A_1,\ldots, A_h$ are sets of positive integers 
with $\sigma(A_i) = \alpha_i$ for all $i \in \{1,\ldots, h\}$, then 
\[
1 - \sigma(A_1+\cdots + A_h)  \leq \prod_{i=1}^h (1-\alpha_i).
\]
If $A_i = A$ for all $i  \in \{1,\ldots, h\}$ and $\sigma(A) = \alpha$, then 
\[
1 - \sigma(hA) \leq (1-\alpha)^h
\]
or 
\[
\sigma(hA) \geq 1 - (1-\alpha)^h
\]
If $\alpha > 0$, then there exists $h_0$ such that 
$\sigma(h_0A) \geq 1/2$.  
From Theorem~\ref{Ndensity:theorem:Shnirel-pigeonhole} we obtain 
$\sigma(2h_0A) = 1$.   This proves 

\bt[Shnirel'man's basis theorem]
Let $A$ be a set of positive integers.  
If $\sigma(A) > 0$, then the set $A$ is a basis for the nonnegative integers. 
\et

It is natural to ask if  Shnirel'man's inequality and Shnirel'man's basis theorem 
can be extended to $n$ dimensions, that is, to sums of subsets of  $\N_0^n$,  
the semigroup of $n$-dimensional nonnegative lattice points, 
for all $n \geq 2$.
In a New York Number Theory Seminar talk ``Shnirel'man density and the Dyson transform''
(cf. Nathanson~\cite{nath26a})  
on September 5, 2024, I said that Shirel'man's inequality and Shnirel'man's basis 
theorem were unsolved problems for lattice points in $\N_0^n$.  
\emph{Mirabile dictu}, the following week I found  on my computer some rough notes on 
Shnirel'man density that I had written 10 years ago 
and that included a proof of Shnirel'man's inequality for lattice points.    
These notes contained no references and I do not know if it was my proof or if 
it was a proof that I had found in a published paper or in a document  floating 
in the internet.  
The proof in these old notes was specifically about sums of sets of lattice points, but 
a study of the argument showed that it depended only on  
order-theoretic properties of the group of lattice points 
and that the theorem could be extended to partially ordered abelian and  nonabelian groups.  
This paper is devoted to proving these group theoretic results and recovering 
the addition theorems for $\N_0^n$. 

The first step is to find analogues of Shnirel'man density that can be 
constructed on arbitrary sets and, in particular, on ordered 
or partially ordered abelian and nonabelian groups. 
We can describe classical Shnirel'man density as follows.  
Let $\mcj$ be the set of all intervals of integers of the form $J = \{1, 2, \ldots, n\}$.  
For every set $A$ of positive integers,  
\[
\sigma(A) = \inf\left\{ \frac{|A\cap J|}{|J|}:J \in \mcj \right\}.
\]
This suggests the following construction of a generalized density 
on the set of subsets of an nonempty set  $X$.  
Let $\mcj$ be a set of nonempty finite subsets of $X$.  
For every subset $A$ of $X$, let 
\beq             \label{Ndensity:def}
\sigma_{\mcj}(A) = \inf \left\{ \frac{|A \cap J|}{|J|} : J \in \mcj \right\}.
\eeq
This density has the following properties.

\bt
Let $X$ be a set and let $\mcj$ be a set of nonempty finite subsets of $X$.  
Let $A$ and $A'$ be subsets of $X$. 
\benum
\item[(i)]
If $A \subseteq A'$, then 
\[
0 \leq \sigma_{\mcj}(A)  \leq \sigma_{\mcj}(A') \leq 1.
\]
\item[(ii)]
$\sigma_{\mcj}(A) = 0$ if $A\cap J = \emptyset$ for some $J \in \mcj$.
\item[(iii)]
$\sigma_{\mcj}(A) = 1$ if and only if $\bigcup_{J \in \mcj} J \subseteq A$.
\item[(iv)]
If $X = \bigcup_{J \in \mcj} J $, then $\sigma_{\mcj}(A) = 1$ if and only if $A = X$.
\eenum
\et

\begin{proof}
If $A \subseteq A'$ and $J \in \mcj$, then $A\cap J \subseteq A' \cap J \subseteq J$ and so 
\[
0 \leq \frac{|A \cap J|}{|J|} \leq \frac{|A' \cap J|}{|J|} \leq 1.
\]
This proves~(i). 

The definition of $\sigma_{\mcj}(A)$ immediately implies~(ii). 

To prove~(iii), it suffices to observe that $\sigma_{\mcj}(A) = 1$ if and only if 
\[
\frac{|A \cap J|}{|J|} =1
\]
for all $J \in \mcj$ if and only if $J\subseteq A$ for all $J \in \mcj$ 
if and only if $\bigcup_{J \in \mcj} J \subseteq A$.   
This implies~(iv) and completes the proof. 
\end{proof}

\section{Partially ordered groups} 
A \emph{partial order} on a nonempty set $X$ is a relation $\leq$ that satisfies 
the following properties: 
\benum
\item
Reflexivity: $x \leq x$ for all $x \in X$, 
\item
Anti-symmetry: For all $x,y \in X$, if $x \leq y$ and $y \leq x$, then $x=y$, 
\item
Transitivity:  For all $x,y,z \in X$, if $x \leq y$ and $y \leq z$, then $x \leq z$.
\eenum
We write $x < y$ if $x \leq y$ and $x \neq y$.
A \emph{total order} on $X$ is a partial order $\leq$ 
such that,  for all $x,y \in X$, either $x \leq y$ or $y \leq x$.

The following is a fundamental result in the theory of ordered sets.

\bt[Szpilrajn~\cite{szpi30}]       \label{lattice:theorem:Szpilrajn}
Every partial order on a set $X$ can be extended to a total order on $X$.
\et

For a short and beautiful proof of Szpilrajn's theorem, see Mandler~\cite{mand20}.

A \emph{lower bound} for a nonempty subset $Y$ of a partially ordered set $X$ 
is an element $y_0 \in X$ such that $y_0 \leq y$ for all $y \in Y$. 
A partially ordered set $X$ is \emph{well-ordered} if every 
nonempty subset of $X$ contains a lower bound.

Let $G$ be a group  with identity element $e$ 
and let $\leq$ be a partial order on the set $G$.  
The following definitions are standard in the theory of partially ordered groups 
(cf.   Clay-Rolfson~\cite{clay-rolf01} and Glass~\cite{glas99}). 

The group $G$ is \emph{right partially ordered} if $a \leq b$ implies $ac \leq bc$ 
for all $a,b,c \in G$.

The group $G$ is \emph{left partially ordered} if $a \leq b$ implies $ca \leq cb$ 
for all $a,b,c \in G$. 

The group $G$ is \emph{partially ordered} if it is both left and right partially ordered. 

In a right or left partially ordered group, if $a > e$ and $b > e$, then $ab > e$. 

The group $G$ is \emph{ordered} if it is partially ordered  and the 
partial order is a total order. 

Let $G$ be an ordered group.  If $x \in G$ and $e < x$, then 
$e < x < x^2 < x^3 < \cdots$ and so $x$ has infinite order.  
Thus, every ordered group is torsion-free, but 
not every torsion-free group is an ordered group.  
However, Levi~\cite{levi42,levi43} proved that  torsion-free abelian groups  
and free groups of finite rank are orderable .  

A subgroup of a right (resp. left) partially ordered group 
is a  right (resp.  left) partially ordered group.  
If $(G_i,\leq_i)$ is a right (resp. left) partially ordered group for all $i \in I$, 
then the direct product $\prod_{i\in I} G_i$ is a right (resp. left) partially ordered group 
with the \emph{rectangular partial order}  
defined by $(a_i)_{i\in I} \leq (b_i)_{i\in I}$ if $a_i \leq_i b_i$ for all $i \in I$. 

Let $X$ be a well-ordered set and let $G$ be a totally ordered additive abelian group. 
The set $G^X$ of functions $f:X \rightarrow G$ is an additive abelian group 
with addition defined by $(f+g)(x) = f(x) + g(x)$. 
The \emph{lexicographical order} on $G^X$ is a total order 
that is defined as follows: 
$f \lex g$ if $f=g$ or if $f \neq g$ and $f(y_0) < g(y_0)$, where 
$y_0 \in Y$ is the lower bound of the nonempty set $Y = \{y\in X: f(y) \neq g(y)\}$. 
The element $y_0 \in Y$ exists because the set $X$ is well-ordered. 

For example, in the additive abelian group $\Z^n$, let 
\[
\mbx = \vectorxn, \qquad  \mby =  \vectoryn, \qquad  \mbz =  \vectorzn.  
\]
Denote the zero vector in $\Z^n$ by $\mbo$. 
The rectangular partial order in $\Z^n$ is defined by $\mbx \leq \mby$ 
if $x_i \leq y_i$ for all $i \in \{1,\ldots, n \}$.
For all $\mbx, \mby \in \Z^n$,  the open interval 
\begin{align*}
( \mbx, \mby) & = \left\{ \mbz \in \Z^n: \mbx < \mbz < \mby \right\} \\ 
& = \left\{ \mbz \in \Z^n \setminus \{\mbx,\mby\} : x_i  \leq z_i  \leq y_i \text{ for all } i \in \{1,\ldots, n\} \right\} 
\end{align*}
is an $n$-dimensional ``rectangle'' of  cardinality $\prod_{i=1}^n (y_i - x_i + 1) -2$. 
Thus, every open interval in $\Z^n$ is finite with respect to the rectangular partial order. 

Every lattice point in $\Z^n$ is a function from the well-ordered set 
$\{1,\ldots, n\}$ into \Z.
The \emph{lexicographical order} on $\Z^n$ is defined as follows:
 $\mbx \lex \mby$ if $\mbx = \mby$ or if 
$x_i < y_i $, where $i$ is the smallest integer 
in $\{1,\ldots, n\}$ such that $x_i \neq y_i$.
The lexicographical  order in $\Z^n$ is a total order 
that extends the rectangular order 
in the sense that $\mbx \leq \mby$ implies $\mbx \lex \mby$ for all 
$\mbx,\mby \in \R^n$.  But not conversely: For the vectors 
$\mbe_1= \bmat 1\\ 0 \emat$ and $\mbe_2 = \bmat 0 \\ 1 \emat$ in $\R^2$, 
we have $\mbe_2 \lex \mbe_1$ but $\mbe_2 \not\leq \mbe_1$. 
The extension of $\leq$ by $\lex$  is an instance of Szpilrajn's theorem.  

Let $G$ be a right or left partially ordered group with identity element $e$.     
For $x,y \in G$, we define the open interval 
\[
(x,y) = \{z \in G: x < z < y \}.
\]
The \emph{positive cone} in $G$ is the set 
\[
G^+ = \{x\in G: x > e\}.
\]
This cone determines the order in $G$. 
If $G$ is a right partially ordered cone, then $a < b$ if and only if $ba^{-1} \in G^+$. 
If $G$ is a left partially ordered cone, then $a < b$ if and only if $a^{-1}b \in G^+$. 

For $x \in G^+$, we define the open interval 
\[
(e,x) = \{y \in G^+: e < y < x \}.
\]
A subset $J$ of $G^+$ is \emph{downward closed} if $x \in J$ implies 
$(e,x) \subseteq J$.

Let $A$ and $B$ be subsets of a group $G$.
We define the \emph{product set} 
\[
AB =  \{ab:a\in A \cup \{e\} \text{ and } b \in B \cup \{e\} \}.
\]
For $h \geq 2$, the product of $h$  subsets $A_1,\ldots, A_h$ of $G$ 
is 
\[
A_1\cdots A_h = \{a_1\cdots a_h: a_i \in A_i \cup \{e\} \text{ for all } i \in \{1,\ldots, h\} \}.
\]
If $A_i = A$ for all $i \in \{1,\ldots,, h\}$, then 
\[
A^h = \underbrace{A\cdots A}_{\text{$h$ factors}} 
= \left\{ a_1\cdots a_h: a_i \in A \cup \{e\} \text{ for all } i \in \{1,\ldots, h\} \right\}.
\]
Let $X$ be a subset of $G$.  
The set $A$ is a \emph{basis of order $h$}  for $X$ if $X \subseteq A^h$.  
The set $A$ is a \emph{basis} for $X$ if $X \subseteq A^h$ 
for some positive integer $h$.

The following results extend Theorem~\ref{Ndensity:theorem:Shnirel-pigeonhole} 
to partially ordered groups. 

\bt                       \label{Ndensity:theorem:J-pigeonhole} 
Let  $G$ be a partially ordered group with positive cone $G^+$. 
Let \mcj\ be a set of  nonempty finite subsets of $G^+$.  
Let $\sigma_{\mcj}$ be the density defined on $G$ by~\eqref{Ndensity:def} and  
let $A$ and $B$ be subsets of $G^+$ with $\sigma_{\mcj}(A) + \sigma_{\mcj}(B)> 1$.  
If $x \in G^+$ and $(e,x) \in \mcj$, then there exist $a \in A$ and $b \in B$ with $ab=x$. 
\et

\begin{proof}
The open interval $(e,x)$ is finite because $(e,x) \in \mcj$.  
Let $J = (e,x)$.  
From the definition of the density $\sigma_{\mcj}$, we have 
\[
|A\cap J| \geq \alpha |J| \qqand |B\cap J| \geq \beta |J|.
\]
If $b \in B \cap J$, then  $e < b < x$.  
Right multiplication of $e < b < x$ by $b^{-1}$ gives $b^{-1} < e < xb^{-1}$. 
Left multiplication of $b^{-1} < e$ by $x$ gives $ xb^{-1} < x$. 
Thus, if $b \in B \cap J$, then $e < xb^{-1} < x$ and so $xb^{-1} \in J$.  
Therefore,    
\[
x(B \cap J)^{-1} \subseteq J    
\]
and so 
\[
\left| x(B \cap J)^{-1}  \right| = \left| B \cap J \right| \geq \beta |J|. 
\]
Both $A\cap J$ and $x(B \cap J)^{-1}$ are subsets of $J$ and   
\[
|A\cap J| + \left| x(B \cap J)^{-1}  \right| \geq (\alpha+\beta)|J| > |J|.   
\]
By the pigeonhole principle,  the sets $A\cap J$ and $x(B \cap J)^{-1}$ 
of $J$ are not disjoint and so there exist $a \in A\cap J$ and $b \in B\cap J$ 
such that $a = xb^{-1}$.    Thus, $x = ab \in AB$. 
This completes the proof. 
\end{proof}

An element $x \in G^+$ is an \emph{atom} if $(e,x) = \emptyset$.  
Equivalently, $x \in G^+$ is an atom if $e < x$ and there exists no $y \in G^+$ 
such that $e < y < x$. 

\bt                           \label{Ndensity:theorem:J-pigeonhole-2} 
Let  $G$ be a partially ordered group with positive cone $G^+$.  
Let $H^+$ be the set of all atoms in $G^+$. 
Let \mcj\ be a set of  nonempty finite subsets of $G^+$ such that 
$(e,x) \in \mcj$ for all $x \in G^+ \setminus H^+$.   
Let $A$ and $B$ be subsets of $G^+$ with $H^+ \subseteq A \cup B$.  
If $\sigma_{\mcj}(A) + \sigma_{\mcj}(B)> 1$, then $AB = G^+ \cup \{e\}$. 
\et

\begin{proof}
We have $H^+ \subseteq A\cup B \subseteq AB$.  
If $x \in G^+\setminus H^+$, then $(e,x) \in \mcj$ and so $x \in AB$ by 
Theorem~\ref{Ndensity:theorem:J-pigeonhole}.  
This completes the proof. 
\end{proof}

\section{Addition theorems for groups}

In this section we prove an analogue of Shnirel'man's inequality (Theorem~\ref{Ndensity:theorem:Shnirel-ineq-classical}) for partially ordered 
abelian and nonabelian groups.

\bt                     \label{Ndensity:theorem:NonabelianPartition}
Let $G$ be a right or left partially ordered group with positive cone 
$G^+ = \{x\in G: x > e\}$. 
Let  $B$ be a nonempty subset of $G^+$. 
Let $J$   be a downward closed  subset of $G^+$ such that 
$J^* = J \setminus B$ is nonempty and, for all $x \in J^*$,  the set  
\[
B^*(x) = \{b \in B:b < x\} = (e,x) \cap B 
\]
is nonempty and finite. 
Then there is a set $\{ J_{\ell} : \ell \in L\}$ of pairwise disjoint nonempty subsets of $J^*$ 
such that 
\[
J^* = \bigcup_{\ell\in L} J_{\ell} 
\]
and, for all $\ell \in L$, there exists $b_{\ell} \in B$ such that 
\benum
\item[(i)]
if $G$ is right partially ordered, then 
$J_{\ell} b_{\ell}^{-1} $ is a downward closed subset of the cone $G^+$,
\item[(ii)]
if $G$ is left partially ordered, then 
$b_{\ell}^{-1} J_{\ell}$ is a downward closed subset of the cone $G^+$. 
\eenum
\et

Note:  A stronger condition on the group $G$ is that every interval $(e,x)$ is finite and 
that $(e,x) \cap B \neq \emptyset$ for all $x \in J^*$.  The condition that open intervals 
be finite is a  discreteness condition on a partially ordered group.  
Here are two examples of groups in which every interval is infinite.  
Consider the additive group \R\ with the usual ordering.  
If $G$ is a subgroup of \R\ that contains \Q, 
then $(0,x)$ is infinite for all $x > 0$.  
If $G$ is an ordered divisible group, then for all $x > e$ and for all $ k \in \N$ 
there exists $y \in G$ such that $x = y^k$ and so 
\[
e < y < y^2 < \cdots < y^k = x
\]
and $|(e,x)| \geq k-1$.  Thus, the open interval $(e,x)$ in $G$ is infinite for all $x > e$.

\begin{proof}
By Szpilrajn's theorem (Theorem~\ref{lattice:theorem:Szpilrajn}), 
there is a total order $\preceq_{tot}$ on $G$ that extends $\leq$. 
Note that we do not assume that $G$ is a right or left partially ordered group 
with respect to $\preceq_{tot}$.

Let $G$ be a  right   partially ordered group with respect to $\leq$ 
and let $x \in J^* = J \setminus B$.  
Because the set  $B^*(x)$ is nonempty and  finite  and 
because $\preceq_{tot}$ is a total order, there exists 
$b^*(x) \in B^*(x)$ such that $b \preceq_{tot} b^*(x)$ for all $b \in B^*(x)$. 
Let 
\[
\{ b^*(x): x \in J^*\} = \{ b_{\ell}:\ell  \in L \} 
\]
where $b_{\ell _1} \neq b_{\ell _2}$ for all $\ell _1,\ell _2 \in L$ 
with $\ell _1 \neq \ell _2$. 
For all $\ell  \in L$, let 
\[
J_{\ell} = \left\{x \in J^*: b^*(x) =  b_{\ell } \right\}. 
\]
For each $x \in J^*$ there is a unique $\ell  \in L$ such that 
$b^*(x) = b_{\ell }$ and so the sets $\{J_{\ell}:\ell  \in L\}$ are pairwise disjoint and 
\[
J^* = \bigcup_{\ell \in L} J_{\ell}.
\]

If $x \in J_{\ell}$, then $b_{\ell} < x$ and so $e < x b_{\ell}^{-1}$.  
Thus, $J_{\ell} b_{\ell}^{-1}$ is a subset of $G^+$. 
To show that  $J_{\ell} b_{\ell}^{-1}$ is a downward closed   
subset of $G^+$, 
we must prove that if $a \in G^+$ and $a <  x b_{\ell}^{-1} $ 
for some $x \in J_{\ell}$, then $a \in J_{\ell} b_{\ell}^{-1}$.
Equivalently, we must prove that $w = a b_{\ell}  \in J_{\ell}$.  

Because $e < a$ and $e < b_{\ell}$, we have 
\[
e < b_{\ell}< a b_{\ell} = w < x.
\]
The set $J$ is a downward closed subset of $G^+$ and 
\[
x \in J_{\ell} \subseteq J^* \subseteq J.  
\] 
It follows that 
\[
w \in J.
\]
Because the total order $\preceq_{tot}$ extends the partial order $\leq$, 
the inequality $ b_{\ell} < w $ also implies $b_{\ell} \prec_{tot}  w $.  
 If $w \in J \cap B$,  then $w < x$ implies $w \in B^*(x)$ and so 
\[
b_{\ell} \prec_{tot} w \preceq_{tot} b^*(x) = b_{\ell}
\]
which is absurd. Therefore,  $w \in J\setminus B = J^*$.  

Because $x \in J_{\ell}$, the inequality 
\[
b_{\ell} < w < x
\]
implies 
\[
b_{\ell} \preceq_{tot}  b^*(w) \preceq_{tot}  b^*(x)= b_{\ell} 
\]
and so 
\[
b_{\ell} = b^*(w) = b^*(x). 
\]
Therefore, $a b_{\ell} = w \in J_{\ell}$ and 
$J_{\ell} b_{\ell}^{-1}$ is a downward closed subset of $G^+$ 
for all $\ell  \in L$. 
This completes the proof when $G$ is right partially ordered.  
The proof is similar when $G$ is left  partially ordered.
\end{proof}

\bt	            \label{Ndensity:theorem:OrderShnirelman}  
Let $G$ be a   right or left partially ordered group with positive cone 
$G^+ = \{x\in G: x > e \}$ and let \mcj\ be the set of all downward closed 
nonempty finite subsets of $G^+$. 
Let $\sigma_{\mcj}$ be the density defined by \mcj\ on subsets of $G^+$. 
Let $A$ and $B$ be subsets of $G^+$ with densities 
$\sigma_{\mcj}(A) = \alpha$ and $\sigma_{\mcj}(B) = \beta$.  
Suppose that, for all $J \in \mcj$ and $x \in J\setminus B$, the set 
\[
B^*(x) = \{b \in B: b < x\} = B \cap (e,x) 
\]
is nonempty and finite.  
\benum
\item[(i)]
If $G$ is right partially ordered, then 
\[ 
\sigma_{\mcj}(AB) \geq \alpha + \beta - \alpha\beta.
\]
\item[(ii)]
If $G$ is left partially ordered, then 
\[ 
\sigma_{\mcj}(BA) \geq \alpha + \beta - \alpha\beta.
\]
\eenum
\et

\begin{proof}
Let $AB = C$.  Recall that $A \subseteq C$ and $B \subseteq C$.  We must prove that
\beq		       \label{Ndensity:ineqCJ}
\frac{|C\cap J|}{|J|} \geq \alpha + \beta - \alpha\beta   
\eeq
for all $J \in \mcj$. 
    
If $J \in \mcj$ and $J \subseteq B$, then $J = B \cap J = C \cap J$ 
and  
\[
\frac{|C\cap J|}{|J|} =  \frac{ |J|}{|J|} =  1 \geq \alpha + \beta - \alpha\beta. 
\]
Thus, we can assume that $J \not\subseteq B$ and so 
$J^* = J \setminus B$
is a nonempty   subset of $G^+$.
By Theorem~\ref{Ndensity:theorem:NonabelianPartition}, 
if $G$ is a right partially ordered group, 
then there is a set 
of pairwise disjoint nonempty finite sets $\left\{J_{\ell} : \ell \in L\right\}$
such that $J^* = \bigcup_{\ell\in L} J_{\ell} $ and, for all $\ell \in L$, 
an element $b_{\ell} \in B$ such that the set $J_{\ell} b_{\ell}^{-1} $ is a downward closed  
nonempty finite  subset of $G^+$, that is, $J_{\ell} b_{\ell}^{-1}  \in \mcj$.

The set $J$ is the disjoint union of the following finite sets: 
\[
J = (B \cap J) \cup J^*
= (B \cap J) \cup  \bigcup_{\ell\in L} J_{\ell}.
\]
Since $B \subseteq C$, we have $B\cap J \subseteq C \cap J$, 
and $C \cap J$ is the disjoint union of the following finite sets:  
\[
C\cap J = (B \cap J)  \cup \bigcup_{\ell \in L}\left(C \cap J_{\ell}\right).
\]
We shall compute a lower bound for the cardinality of $C \cap J_{\ell}$.
If
\[
a \in A \cap J_{\ell}  b_{\ell}^{-1} ,
\]
then
\[
a   b_{\ell} \in J_{\ell} \qqand a   b_{\ell}  \in C
\]
and so
\[
a b_{\ell}  \in C \cap J_{\ell}.
\]
Therefore, 
\[ 
\left( A \cap \left(J_{\ell} b_{\ell}^{-1} \right) \right)b_{\ell} \subseteq C \cap J_{\ell}.
\]
Because $J_{\ell} b_{\ell}^{-1} \in \mcj$ and $\sigma_{\mcj}(A) = \alpha$, we have 
\begin{align*}
|C \cap J_{\ell}|
& \geq  \left|  \left( A \cap \left(J_{\ell} b_{\ell}^{-1} \right) \right)b_{\ell}  \right| \\
& =  \left| A \cap \left( J_{\ell} b_{\ell}^{-1}\right) \right| \\
& \geq  \alpha \left| J_{\ell} b_{\ell}^{-1}\right| \\
& =  \alpha |J_{\ell} |. 
\end{align*}
It follows that  
\begin{align*}
|C\cap J| & =   |B \cap J| + \sum_{\ell \in L} |C \cap J_{\ell}| \\
& \geq   |B \cap J| + \alpha \sum_{\ell \in L} |J_{\ell}| \\
& =   |B \cap J| + \alpha |J^*| \\
& =   |B \cap J| + \alpha \left( |J| - |B \cap J|\right) \\
& =   (1-\alpha)|B \cap J| + \alpha |J| \\
& \geq   (1-\alpha)\beta |J| + \alpha |J| \\
& =   (\alpha + \beta - \alpha\beta) |J|
\end{align*}
and so 
\[
\frac{|C\cap J|}{|J|} \geq \alpha + \beta - \alpha\beta  
\]
and 
\[
\sigma_{\mcj}(C) = \inf_{J\in \mcj} \frac{|C\cap J|}{|J|}
\geq \alpha + \beta - \alpha\beta.
\]
This completes the proof when $G$ is right partially ordered.  
The proof is similar when $G$ is left partially ordered. 
\end{proof}

\bt	            \label{Ndensity:theorem:OrderShnirelman-inductive}  
Let $G$ be a partially ordered group with positive cone 
$G^+ = \{x\in G: x > e \}$ and let \mcj\ be the set of all downward closed 
nonempty finite subsets of $G^+$.   
Let $\sigma_{\mcj}$ be the density defined by \mcj\ on subsets of $G^+$. 
For $h \geq 2$, let $A_1, \ldots, A_h$ be subsets of $G^+$ with  
$\sigma_{\mcj}(A_i) = \alpha_i$ for all $i \in \{1,\ldots, h\}$.    
Suppose that, for all $J \in \mcj$, for all $i \in \{1,\ldots, h\}$, 
and for all $x \in J\setminus A_i$, the set 
\[
A_i^*(x) = \{a_i \in A_i: a_i < x\} = A_i \cap (e,x) 
\]
is nonempty and finite.  
Then    
\[
1 - \sigma_{\mcj}\left(A_1\cdots A_h \right) 
\leq \prod_{i=1}^h \left(1 -  \sigma_{\mcj}\left( A_i \right) \right). 
\]
\et

\begin{proof}
This follows from Theorem~\ref{Ndensity:theorem:OrderShnirelman} 
by induction on $h$. 
\end{proof}

\bt	            \label{Ndensity:theorem:h-OrderShnirelman-basis}  
Let $G$ be a partially ordered group,  
let $G^+ = \{x\in G: x > e \}$ be the positive cone in $G$, and   
let \mcj\ be the set of all downward closed nonempty finite subsets of $G^+$.
Let $A$ be a nonempty subset of $G^+$ such that, 
for all $J \in \mcj$ and $x \in J\setminus A$, the set 
\[
A^*(x) = \{a \in A: a < x\} = A \cap (e,x) 
\]   
is nonempty and finite.  
Let $\sigma_{\mcj}$ be the density defined by \mcj\ on subsets of $G^+$. 
If  $\sigma_{\mcj}(A) = \alpha > 0$, then $A$ is a basis for $G$. 
\et

\begin{proof}
Let $x \in G^+$.  If $x$ is an atom in $G$, that is, 
if the open interval $(e,x)$ is empty, then the set $\{x\}$ is 
a downward closed nonempty finite subset of $G^+$ and so $\{x\} \in J$.  
We have 
\[
0 < \sigma_{\mcj}(A) \leq \frac{|A \cap \{x\}|}{| \{x\}|} = |A \cap \{x\}| 
\]
and so $x \in A \subseteq A^h$ for all $h \geq 2$. 

Let $h \geq 2$.  Applying Theorem~\ref{Ndensity:theorem:OrderShnirelman-inductive} 
with $A_i = A$ 
for all $i \in \{1,\ldots, h\}$, we obtain 
\[
1 - \sigma_{\mcj}(A^h) \leq (1 - \alpha)^h 
\]
and so there exists $h_0$ such 
\[
\sigma_{\mcj}(A^{h_0}) > \frac{1}{2}. 
\]
If $(e,x)$ is not empty, then, by Theorem~\ref{Ndensity:theorem:J-pigeonhole-2}, 
there exist $a,b \in A^{h_0}$ such that $x=ab \in A^{2h_0}$ 
and so $A$ is a basis for $G$ of order $2h_0$.    
This completes the proof. 
\end{proof}

\section{Sums of sets of lattice points}
We return to additive number theory for lattice points.  
Let $\Z^n$ be the partially ordered additive abelian group of lattice points 
with the rectangular order.  
The positive cone of $\Z^n$ is the semigroup  $\N_0^n\setminus \{\mbo\}$ 
of nonzero nonnegative lattice points.  

The set $A$ of nonnegative lattice points is a 
\emph{basis of order $h$} for $\N_0^n \setminus \{\mbo\}$ 
if $hA = \N_0^n \setminus \{\mbo\}$.  
The set $A$ is a \emph{basis} for $\N_0^n \setminus \{\mbo\}$ 
if $A$ is a basis of order $h$ for $\N_0^n \setminus \{\mbo\}$ for some $h$.

\bt	            \label{Ndensity:theorem:LatticeShnirelman} 
Let $\Z^n$ be the partially ordered additive abelian group of lattice points 
with the rectangular order and with positive cone $\N_0^n\setminus \{\mbo\}$.   
Let \mcj\ be the set of all downward closed 
nonempty finite subsets of  $\N_0^n\setminus \{\mbo\}$  
and let $\sigma_{\mcj}$ be the density on  $\Z^n$ defined by \mcj. 
If $A$ and $B$ are subsets of $\N_0^n\setminus \{\mbo\}$ with densities 
$\sigma_{\mcj}(A) = \alpha$ and $\sigma_{\mcj}(B) = \beta$, then 
\beq                                                     \label{Ndensity:LatticeShnirelman} 
\sigma_{\mcj}(A + B) \geq \alpha + \beta - \alpha\beta.
\eeq
\et

\begin{proof}
We must prove that 
 \[
 \frac{ \left| (A+B) \cap J \right|}{\left| J \right|} \geq \alpha + \beta - \alpha\beta
\]
for all $J \in \mcj$.
Because $A \subseteq A+B$, if $\beta = 0$, then 
\[
\sigma_{\mcj}(A+B) \geq \sigma_{\mcj}(A) = \alpha = \alpha + \beta - \alpha\beta.
\]
Thus, we can assume that $\beta > 0$.

For all $i \in \{1,\ldots, n\}$, let $\mbe_i$  be the standard unit vector whose $j$th coordinate 
is the Kronecker delta $\delta_{i,j}$.  If $\mbx$ is a nonzero vector in $\N_0^n \setminus \{\mbo\}$, 
then $x_i \geq 1$ for some $i \in \{1,\ldots, n\}$ and so $\mbe_i \leq \mbx$.  
It follows that if $J$ is a downward closed nonempty subset of $\N_0^n \setminus \{\mbo\}$, 
then 
\[
 \{\mbe_1,\ldots, \mbe_n\}  \cap J\neq \emptyset. 
\]

For all $i \in \{1,\ldots, n\}$, the unit vector $\mbe_i$ is an atom 
and the set $\{\mbe_i\}$ is a downward closed finite nonempty 
subset of $\N_0^n \setminus \{\mbo\}$.  Thus, $\{\mbe_i\} \in \mcj$.   
We have 
\[
0 <\beta = \sigma_{\mcj}(B) \leq \frac{ \left| B \cap \{\mbe_i\} \right|}{\left|  \{\mbe_i\} \right|}
= \begin{cases}
1 & \text{if $\mbe_i \in B$} \\ 
0 & \text{if $\mbe_i \notin B$ }
\end{cases}
\]
and so $\mbe_i \in B$ and 
\[
\{\mbe_1,\ldots, \mbe_n\} \subseteq B.  
\]
Therefore, 
 $B \cap J \neq \emptyset$ for all $J \in \mcj$. 
 If $B \cap J = J$, then $B \subseteq A+B$ implies 
 \[
 \frac{ \left| (A+B) \cap J \right|}{\left| J \right|} \geq  \frac{ \left| B \cap J \right|}{\left| J \right|} 
 = \frac{ \left| J \right|}{\left| J \right|} = 1 \geq \alpha + \beta - \alpha\beta.
\]
Thus, we can assume that $B \cap J \neq J$ and so $J \setminus B \neq \emptyset$.  
If $\mbx \in J\setminus B$, then  $\mbx > \mbe_i > \mbo$ 
for some $i \in \{1,\ldots, n\}$ and  $\mbe_i \in B^*(\mbx) = B  \cap (\mbo,\mbx)$.
Every interval in $\Z^n$ is finite and so $B^*(\mbx)$ is nonempty and finite 
for all $\mbx \in J\setminus B$. 
It follows from Theorem~\ref{Ndensity:theorem:OrderShnirelman} that  
$\sigma_{\mcj}(A + B) \geq \alpha + \beta - \alpha\beta.$.  
This completes the proof. 
\end{proof}

\bt	            \label{Ndensity:theorem:LatticeBasis} 
Let $\Z^n$ be the partially ordered additive abelian group of lattice points 
with the rectangular order.  
Let \mcj\ be the set of all downward closed 
nonempty finite subsets of $\N_0^n\setminus \{\mbo\}$  
and let $\sigma_{\mcj}$ be the density on  $\Z^n$ defined by \mcj. 
If $A$ is a  subset of $\N_0^n\setminus \{\mbo\}$ with density 
$\sigma_{\mcj}(A) > 0$, then $A$ is a   basis for 
the semigroup $\N_0^n\setminus \{\mbo\}$. 
\et

\begin{proof}
This follows from Theorems~\ref{Ndensity:theorem:h-OrderShnirelman-basis} 
and~\ref{Ndensity:theorem:LatticeShnirelman}.  
\end{proof}

\section{Open problems}
It is natural to ask what other results about the Shnirel'man density of 
sets of positive integers extend to lattice points and to arbitrary 
abelian and nonabelian partially ordered groups.  
For example, Mann proved an addition theorem for finite sets of integers 
that leads to an important and best possible 
strengthening of Shnirel'man's addition theorem 
(Theorem~\ref{Ndensity:theorem:Shnirel-ineq-classical}).

\bt[Mann~\cite{mann42}]
Let $A$ and $B$ be sets of positive integers with Shnirel'man densities 
$\sigma(A) = \alpha$ and $\sigma(B) = \beta$.  Then 
\[
\sigma(A+B) \geq \min(1,\alpha+\beta).
\]
\et

Dyson~\cite{dyso45} refined Mann's combinatorial result 
about finite sumsets (see Nathanson~\cite{nath26a}).  
Do these  results generalize to partially ordered groups?

\end{document}